\documentclass[12pt,reqno]{amsart}
\usepackage{amsmath,amssymb,amsfonts,amscd,latexsym,amsthm,mathrsfs,verbatim,comment}
\usepackage{hyperref}
\newtheorem{theorem}{Theorem}
\begin{document}

\title{Ramanujan-type formulae for $1/\pi$: $q$-analogues}

\date{18 February 2018}

\author{Victor J. W. Guo}
\address{School of Mathematical Sciences, Huaiyin Normal University, Huai'an, Jiangsu 223300, People's Republic of China}
\email{jwguo@hytc.edu.cn}

\author{Wadim Zudilin}
\address{Department of Mathematics, IMAPP, Radboud University, PO Box 9010, 6500~GL Nijmegen, Netherlands}
\email{w.zudilin@math.ru.nl}

\address{School of Mathematical and Physical Sciences, The University of Newcastle, Callag\-han, NSW 2308, Australia}
\email{wadim.zudilin@newcastle.edu.au}

\thanks{The first author was partially supported by the National Natural Science Foundation of China (grant 11771175).}


\begin{abstract}
The hypergeometric formulae designed by Ramanujan more than a century ago for efficient approximation of $\pi$, Archimedes' constant,
remain an attractive object of arithmetic study. In this note we discuss some $q$-analogues of Ramanujan-type evaluations and of related supercongruences.
\end{abstract}

\maketitle

\section{Introduction}
\label{s1}

Let $q$ be inside the unit disc, $|q|<1$.
In the recent joint paper \cite{GL18} of one of these authors, Ramanujan's formulas \cite{Ra14}
\begin{equation}
\sum_{n=0}^\infty (6n+1)\frac{(\frac{1}{2})_n^3}{n!^3 4^n} =\frac{4}{\pi}
\quad\text{and}\quad
\sum_{n=0}^\infty (-1)^n(6n+1)\frac{(\frac{1}{2})_n^3}{n!^3 8^n } =\frac{2\sqrt{2}}{\pi}
\label{a0}
\end{equation}
were supplied with $q$-analogues
\begin{align}
\sum_{n=0}^{\infty}q^{n^2}[6n+1]\frac{(q;q^2)_n^2(q^2;q^4)_n}{(q^4;q^4)_n^3}
&=\frac{(1+q)(q^2;q^4)_{\infty}(q^6;q^4)_{\infty}}{(q^4;q^4)_{\infty}^2},
\label{a1}
\\
\sum_{n=0}^{\infty}(-1)^n q^{3n^2}[6n+1]\frac{(q;q^2)_n^3}{(q^4;q^4)_n^3}
&=\frac{(q^3;q^4)_\infty (q^5;q^4)_\infty}{(q^4;q^4)_\infty^2}.
\label{a11}
\end{align}
Here and in what follows
$$
(a;q)_\infty=\prod_{j=0}^\infty(1-aq^j)
$$
and we use the standard notation
$$
(a)_n=\frac{\Gamma(a+n)}{\Gamma(a)}
\quad\text{and}\quad
(a;q)_n=\frac{(a;q)_\infty}{(aq^n;q)_\infty}
$$
for the Pochhammer symbol and its $q$-version, so that
$$
(a)_n=\prod_{j=0}^{n-1}(a+j)
\quad\text{and}\quad
(a;q)_n=\prod_{j=0}^{n-1}(1-aq^j)
$$
for \emph{positive} integers $n$. We also define the $q$-numbers by $[n]=[n]_q=(1-q^n)/(1-q)$.

Apart from another entry
$$
\sum_{n=0}^{\infty}(-1)^nq^{n^2}[4n+1]\frac{(q;q^2)_n^3}{(q^2;q^2)_n^3}
=\frac{(q;q^2)_\infty(q^3;q^2)_\infty}{(q;q^2)_\infty^2},
$$
which follows trivially from a limiting case of Jackson's formula \cite[eq.~(2.7.1)]{GR04}
and which represents a $q$-analogue of Bauer's (Ramanujan-type) formula
$$
\sum_{n=0}^\infty (-1)^n(4n+1)\frac{(\frac{1}{2})_n^3}{n!^3} =\frac{2}{\pi},
$$
the entries \eqref{a1} and \eqref{a11} provide us with the first examples of $q$-analogues of Ramanujan's formulae for $1/\pi$.

The principal ingredients in the proof of \eqref{a1} and \eqref{a11} in \cite{GL18} are suitable chosen $q$-Wilf--Zeilberger (WZ) pairs
and some basic hypergeometric identities.

In this note we use the $q$-WZ pairs from \cite{GL18} and original ideas of J.~Guillera from~\cite{Gu06} to give a few further  $q$-analogues of Ramanujan's formulae and their generalizations.

\begin{theorem}
\label{th1}
The following identities are true:
\begin{align}
&
\sum_{n=0}^\infty \frac{(-1)^n q^{2n^2}(q^2;q^4)_{n}^2(q;q^2)_{2n}}{(q^4;q^4)_{n}^2(q^4;q^4)_{2n}}
\biggl([8n+1]+[4n+1]\frac{q^{4n+1}}{1+q^{4n+2}}\biggr)
\nonumber\\ &\qquad
=\frac{(1+q)(q^2;q^4)_{\infty}(q^6;q^4)_{\infty}}{(q^4;q^4)_{\infty}^2},
\label{q1}
\\
&
\sum_{n=0}^\infty \frac{(-1)^n q^{2n^2}(q^2;q^4)_{n}(q^2;q^4)_{2n}(q;q^2)_{3n}}{(q^4;q^4)_{n}^2(q^4;q^4)_{3n}(q;q^2)_n}
\nonumber\\ &\;\times
\biggl([10n+1]+\frac{q^{6n+1}[4n+2][6n+1]}{[12n+4]}+\frac{q^{6n+3}[6n+1][6n+3][8n+2](1+q^{2n+1})}{[12n+4][12n+8]}\biggr)
\nonumber\\ &\qquad
=\frac{(1+q)(q^2;q^4)_{\infty}(q^6;q^4)_{\infty}}{(q^4;q^4)_{\infty}^2},
\label{q2}
\\
&
\sum_{n=0}^\infty \frac{q^{4n^2}(q;q^2)_{2n}^2}{(q^4;q^4)_{n}^2(q^4;q^4)_{2n}}
\biggl([8n+1]-q^{8n+3}\frac{[4n+1]^2}{[8n+4]}\biggr)
=\frac{(q^3;q^4)_\infty (q^5;q^4)_\infty}{(q^4;q^4)_\infty^2}.
\label{q3}
\end{align}
\end{theorem}

As pointed out to us by C.~Krattenthaler,
the formulae \eqref{a1} and \eqref{a11} can be alternatively proved through an intelligent use
of quadratic transformations from \cite[Section 3.8]{GR04}; more details of this machinery appear in \cite{HKS18}.
We adapt the related technique here to establish a new $q$-analogue of Ramanujan's formula as follows.

\begin{theorem}
\label{th2}
We have
\begin{equation}
\sum_{n=0}^{\infty}\frac{q^{2n^2}(q;q^2)_n^2 (q;q^2)_{2n}}{(q^2;q^2)_{2n}(q^6;q^6)_n^2 }[8n+1]
=\frac{(q^3;q^2)_\infty (q^3;q^6)_\infty }{(q^2;q^2)_\infty (q^6;q^6)_\infty}.
\label{q4}
\end{equation}
\end{theorem}

Observe that equations \eqref{q1} and \eqref{q4} are $q$-analogues of Ramanujan's
\begin{gather}
\sum_{n=0}^\infty\frac{(-1)^n(\frac14)_n(\frac12)_n(\frac34)_n}{n!^3\,4^n}\,(20n+3)
=\sum_{n=0}^\infty\frac{(-1)^n{4n\choose 2n}{2n\choose n}^2}{2^{10n}}\,(20n+3)
=\frac{8}{\pi},
\nonumber\\
\sum_{n=0}^\infty\frac{(\frac14)_n(\frac12)_n(\frac34)_n}{n!^3\,9^n}\,(8n+1)
=\sum_{n=0}^\infty\frac{{4n\choose 2n}{2n\choose n}^2}{2^{8n}3^{2n}}\,(8n+1)
=\frac{2\sqrt{3}}{\pi},
\label{non-q4}
\end{gather}
while the $q\to1$ cases of \eqref{q2} and \eqref{q3} read
\begin{gather*}
\sum_{n=0}^\infty \frac{(-1)^n{6n\choose 3n}{4n\choose 2n}{2n\choose n}}{2^{12n}}\frac{576n^3+624n^2+190n+15}{(3n+1)(3n+2)}
=\frac{16}{\pi},
\\
\sum_{n=0}^\infty\frac{{4n\choose 2n}^2{2n\choose n}}{2^{12n}}\frac{48n^2+32n+3}{2n+1}
=\frac{8\sqrt{2}}{\pi}.
\end{gather*}

The details of our proofs of Theorems~\ref{th1} and \ref{th2} are given in Sections~\ref{s2} and \ref{s3}, respectively.
In the final section we highlight some connections of our findings with $q$-analogues of Ramanujan-type supercongruences,
which were an original source of the proofs of \eqref{a1} and \eqref{a11} in~\cite{GL18}.
A simple look of (some) $q$-analogues of Ramanujan-type formulae for $1/\pi$ and certain similarity with the Rogers--Ramanujan identities makes the former a plausible candidate for combinatorial explorations.

As hypergeometric summation and transformation formulae occasionally possess several (sometimes very different!)
$q$-analogues, we cannot exclude a possibility that there are multiple $q$-analogues of some Ramanujan-type formulae for $1/\pi$.
Another source of such multiple $q$-entries can be caused by using different $q$-WZ pairs:
the phenomenon of their existence has been recorded in the non-$q$-settings in Guillera's PhD thesis;
see \cite[Sections~1.3, 1.4]{Gu07}.

\medskip
We thank Jes\'us Guillera for drawing our attention to some parts of his PhD thesis \cite{Gu07} and Christian Krattenthaler for giving us details of his hypergeometric proof of~\eqref{a1}.

\section{WZ machinery and Guillera's invention}
\label{s2}

The heart of the proof of \eqref{a1} in \cite{GL18} is the $q$-WZ pair
\begin{align*}
F(n,k)=F(n,k;q)
&=\frac{q^{(n-k)^2}[6n-2k+1](q^2;q^4)_{n}(q;q^2)_{n-k}(q;q^2)_{n+k}}{(q^4;q^4)_{n}^2(q^4;q^4)_{n-k}(q^2;q^4)_k},
\\
G(n,k)=G(n,k;q)
&=\frac{q^{(n-k)^2}(q^2;q^4)_n (q;q^2)_{n-k}(q;q^2)_{n+k-1}}{(1-q)(q^4;q^4)_{n-1}^2(q^4;q^4)_{n-k}(q^2;q^4)_k},
\end{align*}
where it is set that $1/(q^4;q^4)_{m}=0$ for any negative integer $m$, which satisfies
\begin{align}
F(n,k-1)-F(n,k)=G(n+1,k)-G(n,k).
\label{eq:fnk-gnk}
\end{align}
Summing the both sides of \eqref{eq:fnk-gnk} over $n=0,1,\dots,m-1$ and, afterwards, over $k=1,\dots,m-1$ we get
$$
\sum_{n=0}^{m-1}F(n,0)
=\sum_{k=1}^{m}G(m,k),
$$
equivalently,
\begin{align*}
\sum_{n=0}^{m-1}q^{n^2}[6n+1]\frac{(q;q^2)_n^2(q^2;q^4)_n}{(q^4;q^4)_n^3}
&=\sum_{k=1}^{m}\frac{q^{(m-k)^2}(q^2;q^4)_m(q;q^2)_{m-k}(q;q^2)_{m+k-1}}
{(1-q)(q^4;q^4)_{m-1}^2(q^4;q^4)_{m-k}(q^2;q^4)_{k}}
\\ \intertext{(changing the summation index from $k$ to $m-k$)}
&=\frac{(q^2;q^4)_m}{ (1-q)(q^4;q^4)_{m-1}^2}\sum_{k=0}^{m-1}\frac{q^{k^2}(q;q^2)_k(q;q^2)_{2m-k-1}}
{(q^4;q^4)_k(q^2;q^4)_{m-k}}.
\end{align*}
Now letting $m\to\infty$ and using
$$
\sum_{k=0}^{\infty}q^{k^2}\frac{(q;q^2)_k}{(q^4;q^4)_k}
=\frac{(q^2;q^4)_{\infty}^2}{(q;q^2)_{\infty}},
$$
due to Slater, we arrive at \eqref{a1}.

\medskip
If we denote $\tilde F(n,k)=F(n,-k)$ and $\tilde G(n,k)=G(n,-k)$ then relation \eqref{eq:fnk-gnk} becomes a `standard' form of the WZ relation
\begin{align}
\tilde F(n,k+1)-\tilde F(n,k)=\tilde G(n+1,k)-\tilde G(n,k).
\label{eq:WZ}
\end{align}
It was pointed out by Guillera in \cite{Gu06} that we can iterate a WZ pair $(\tilde F_1,\tilde G_1)$ satisfying \eqref{eq:WZ} by taking
$$
\tilde F_2(n,k)=\tilde F_1(n,n+k)+\tilde G_1(n+1,n+k)
\quad\text{and}\quad
\tilde G_2(n,k)=\tilde G_1(n,n+k).
$$
Indeed,
\begin{align*}
\tilde F_2(n,k+1)-\tilde F_2(n,k)
&=\bigl(\tilde F_1(n,n+k+1)+\tilde G_1(n+1,n+k+1)\bigr)
\\ &\quad
-\bigl(\tilde F_1(n,n+k)+\tilde G_1(n+1,n+k)\bigr)
\\
&=\tilde G_1(n+1,n+k)-\tilde G_1(n,n+k)
\\ &\quad
+\tilde G_1(n+1,n+k+1)-\tilde G_1(n+1,n+k)
\\
&=\tilde G_2(n+1,k)-\tilde G_2(n,k).
\end{align*}
He also observes in \cite{Gu06} that the sums
\begin{equation*}
\sum_{n=0}^\infty\tilde F_1(n,k)=C
\end{equation*}
do not depend on $k=0,1,2,\dots$ (and, in most of the cases, even on $k>-1$ \emph{real} thanks to Carlson's theorem\,---\,this
allows one to compute $C$ by choosing an appropriate real value for~$k$) and that
\begin{equation*}
\sum_{n=0}^\infty\tilde F_2(n,k)=C
\end{equation*}
as well, the same right-hand side.

\begin{proof}[Proof of \eqref{q1} and \eqref{q2}]
The WZ pair $(F,G)$ satisfying \eqref{eq:fnk-gnk} transforms into
\begin{align*}
\tilde F_1(n,k)
&=\frac{(-1)^k q^{n^2+2nk-k^2}[6n+2k+1](q^2;q^4)_{n}(q;q^2)_{n+k}(q;q^2)_{n-k}(q^2;q^4)_{k}}
{(q^4;q^4)_{n}^2(q^4;q^4)_{n+k}},
\\
\tilde G_1(n,k)
&=\frac{(-1)^k q^{n^2+2nk-k^2}(q^2;q^4)_n (q;q^2)_{n+k}(q;q^2)_{n-k-1}(q^2;q^4)_{k}}
{(1-q)(q^4;q^4)_{n-1}^2(q^4;q^4)_{n+k}},
\end{align*}
whose subsequent iterations are
\begin{align*}
\tilde F_2(n,k)
&=\frac{(-1)^n q^{2n^2}(q^2;q^4)_{n}(q;q^2)_{2n+k}(q^2;q^4)_{n+k}}{(q^4;q^4)_{n}^2(q^4;q^4)_{2n+k}(q;q^2)_{k}}\\
&\quad\times\biggl([8n+2k+1]+\frac{q^{4n+2k+1}(1-q^{4n+2})(1-q^{4n+2k+1})}{(1-q)(1-q^{8n+4k+4})}\biggr),
\\
\tilde G_2(n,k)
&=\frac{(-1)^{n-1}q^{2n^2+2k+1}(q^2;q^4)_n (q;q^2)_{2n+k}(q^2;q^4)_{n+k}}
{(1-q)(q^4;q^4)_{n-1}^2(q^4;q^4)_{2n+k}(q;q^2)_{k+1}},
\\ \intertext{and}
\tilde F_3(n,k)
&=\frac{(-1)^n q^{2n^2}(q^2;q^4)_{n}(q;q^2)_{3n+k}(q^2;q^4)_{2n+k}}{(q^4;q^4)_{n}^2(q^4;q^4)_{3n+k}(q;q^2)_{n+k}}
\\
&\quad\times\biggl([10n+2k+1]+\frac{q^{6n+2k+1}(1-q^{4n+2})(1-q^{6n+2k+1})}{(1-q)(1-q^{12n+4k+4})}
\\ &\qquad
+\frac{q^{6n+2k+3}(1-q^{4n+2})(1-q^{6n+2k+1})(1-q^{6n+2k+3})(1-q^{8n+4k+2})}{(1-q)(1-q^{12n+4k+4})(1-q^{12n+4k+8})(1-q^{2n+2k+1})}\biggr),
\\
\tilde G_3(n,k)
&=\frac{(-1)^{n-1}q^{2n^2+2n+2k+1}(q^2;q^4)_n (q;q^2)_{3n+k}(q^2;q^4)_{2n+k}}{(1-q)(q^4;q^4)_{n-1}^2(q^4;q^4)_{3n+k}(q;q^2)_{n+k+1}}.
\end{align*}
It follows \eqref{a1} that in this case
$$
C=\sum_{n=0}^\infty\tilde F_1(n,k)=\sum_{n=0}^\infty\tilde F_1(n,0)
=\frac{(1+q)(q^2;q^4)_{\infty}(q^6;q^4)_{\infty}}{(q^4;q^4)_{\infty}^2},
$$
hence also
$$
\sum_{n=0}^\infty\tilde F_2(n,0)=C
\quad\text{and}\quad
\sum_{n=0}^\infty\tilde F_3(n,0)=C,
$$
which are precisely formulae \eqref{q1} and \eqref{q2}.
\end{proof}

\begin{proof}[Proof of \eqref{q3}]
The main part in proving \eqref{a11} is a WZ pair, whose `standard' form is
\begin{align*}
\tilde F_1(n,k)
=\tilde F_1(n,k;q)
&=(-1)^{n+k}\frac{[6n+2k+1](q;q^2)_{n-k}(q;q^2)_{n+k}^2}{(q^4;q^4)_{n}^2(q^4;q^4)_{n+k}},
\\
\tilde G_1(n,k)
=\tilde G_1(n,k;q)
&=\frac{(-1)^{n+k}(q;q^2)_{n-k-1}(q;q^2)_{n+k}^2}{(1-q)(q^4;q^4)_{n-1}^2(q^4;q^4)_{n+k}},
\end{align*}
with further iteration
\begin{align*}
\tilde F_2(n,k;q)
&=\frac{q^{k^2}(q;q^2)_{2n+k}^2}{(q^4;q^4)_{n}^2(q^4;q^4)_{2n+k}(q;q^2)_k}
\biggl([8n+2k+1]-\frac{[4n+2k+1]^2}{[8n+4k+4]}\biggr),
\\
\tilde G_2(n,k;q)
&=-\frac{q^{(k+1)^2}(q;q^2)_{2n+k}^2}{(1-q)(q^4;q^4)_{n-1}^2(q^4;q^4)_{2n+k}(q;q^2)_{k+1}}.
\end{align*}
The identity \eqref{a11} reads
$$
C=\sum_{n=0}^\infty\tilde F_1(n,0;q^{-1})
=\frac{(q^3;q^4)_\infty (q^5;q^4)_\infty}{(q^4;q^4)_\infty^2}
$$
and leads to
$$
\sum_{n=0}^{\infty}\tilde F_2(n,0;q^{-1})
=C,
$$
which is the desired relation \eqref{q3}.
\end{proof}

\section{A quadratic transformation and a cubic transformation}
\label{s3}

This is the shortest section of this note.

\begin{proof}[New proof of \eqref{a1} and \eqref{a11}]
The formula \cite[eq.~(4.6)]{Ra93} reads
\begin{align}
&
\sum_{n=0}^\infty\frac{(a;q)_n (1-aq^{3n})(d;q)_n(q/d;q)_n(b;q^2)_n}
{(q^2;q^2)_n (1-a)(aq^2/d;q^2)_n (adq;q^2)_n (aq/b;q)_n}\frac{a^n q^{n+1\choose 2}}{b^n}
\notag\\ &\qquad
=\frac{(aq;q^2)_\infty (aq^2;q^2)_\infty (adq/b;q^2)_\infty (aq^2/bd;q^2)_\infty}
{(aq/b;q^2)_\infty (aq^2/b;q^2)_\infty (aq^2/d;q^2)_\infty (adq;q^2)_\infty}.  \label{quadratic}
\end{align}
Letting $q\to q^2$ and taking $a=d=q$ and $b=q^2$, we are led to \eqref{a1}.
Similarly, letting $q\to q^2$ and $b\to\infty$, and taking $a=d=q$, we conclude with \eqref{a11}.
\end{proof}

\begin{proof}[Proof of \eqref{q4}]
Letting $d\to0$ in the cubic transformation \cite[eq.~(3.8.18)]{GR04} we arrive at the formula
\begin{align*}
&
\sum_{n=0}^\infty\frac{(1-acq^{4n})(a;q)_n(q/a;q)_n (ac;q)_{2n}}
{(1-ac)(cq^3;q^3)_n(a^2cq^2;q^3)_n(q;q)_{2n}}\,q^{n^2}
\\ &\qquad
=\frac{(acq^2;q^3)_\infty(acq^3;q^3)_\infty(aq;q^3)_\infty(q^2/a;q^3)_\infty}
{(q;q^3)_\infty(q^2;q^3)_\infty(a^2cq^2;q^3)_\infty(cq^3;q^3)_\infty}.
\end{align*}
Now replace $q$ with $q^2$, and take $a=q$ and $c=1$. The resulting identity is \eqref{q4}.
\end{proof}

\section{$q$-Supercongruences}
\label{s4}

In this section we briefly highlight some links of our results with the supercongruences of `Ramanujan type'.
The general pattern for them in \cite{Zu09} predicts, for example, that the congruence counterpart of Ramanujan's formula \eqref{non-q4} is
$$
\sum_{k=0}^{p-1}\frac{(\frac14)_k(\frac12)_k(\frac34)_k}{k!^3\,9^k}\,(8k+1)
\equiv p\biggl(\frac{-3}p\biggr)\pmod{p^3}
\quad\text{for $p>3$ prime},
$$
where the Jacobi--Kronecker symbol $\bigl(\frac{-3}{\cdot}\bigr)$ `replaces' the square root of $3$. For this particular entry one also observes experimentally that
$$
\sum_{k=0}^{(p-1)/2}\frac{(\frac14)_k(\frac12)_k(\frac34)_k}{k!^3\,9^k}\,(8k+1)
\equiv p\biggl(\frac{-3}p\biggr)\pmod{p^3}
\quad\text{for $p>3$ prime},
$$
when the sum on the left-hand side is shorter.

The $q$-WZ pairs used in the proofs of Section~\ref{s2} were originally designed to verify $q$-analogues of the
supercongruences corresponding to the identities in \eqref{a0}. As suggested by \cite{Zu09},
one can turn the argument into the opposite direction, to provide $q$-analogues of Ramanujan-type supercongruences from Ramanujan-type identities.
Our experimental observations include the $q$-analogues
\begin{align*}
\sum_{k=0}^{(n-1)/2}\frac{(q;q^2)_k^2 (q;q^2)_{2k}}{(q^2;q^2)_{2k}(q^6;q^6)_k^2 }[8k+1]q^{2k^2}
&\equiv q^{-(n-1)/2}[n]\biggl(\frac{-3}{n}\biggr) \pmod{[n]\Phi_n(q)^2},
\\
\sum_{k=0}^{n-1}\frac{(q;q^2)_k^2 (q;q^2)_{2k}}{(q^2;q^2)_{2k}(q^6;q^6)_k^2 }[8k+1]q^{2k^2}
&\equiv q^{-(n-1)/2}[n]\biggl(\frac{-3}{n}\biggr) \pmod{[n]\Phi_n(q)^2}
\end{align*}
for positive $n$ coprime with 6, of the two congruences above, where $\Phi_n(q)$ denotes the $n$-th cyclotomic polynomial.

Similar congruences, in a weaker form, seem to occur for related identities not of Ramanujan type.
For example, replacing $q$ with $q^2$ in \eqref{quadratic}, then choosing $a=q$ and $d=-q$, and 
finally taking $b=q^2$ or $b\to\infty$, respectively, we obtain
\begin{align*}
\sum_{n=0}^\infty\frac{(q^2;q^4)_n (-q;q^2)_n^2}{(q^4;q^4)_n(-q^4;q^4)_n^2}[6n+1]q^{n^2}
&=\frac{(-q^2;q^4)_\infty^2}{(1-q)(-q^4;q^4)_\infty^2},
\\
\sum_{n=0}^\infty(-1)^n\frac{(q;q^2)_n (-q;q^2)_n^2}{(q^4;q^4)_n(-q^4;q^4)_n^2}[6n+1]q^{3n^2}
&=\frac{(q^3;q^4)_\infty (q^5;q^4)_\infty}{(-q^4;q^4)_\infty^2}.
\end{align*}
Numerical experiment suggests the following congruences for the truncated sums:
\begin{align*}
\sum_{k=0}^{(n-1)/2}\frac{(q^2;q^4)_k (-q;q^2)_k^2}{(q^4;q^4)_k(-q^4;q^4)_k^2}[6k+1]q^{k^2}
&\equiv 0\pmod{[n]},
\displaybreak[2]\\
\sum_{k=0}^{n-1}\frac{(q^2;q^4)_k (-q;q^2)_k^2}{(q^4;q^4)_k(-q^4;q^4)_k^2}[6k+1]q^{k^2}
&\equiv 0\pmod{[n]},
\displaybreak[2]\\
\sum_{k=0}^{(n-1)/2}(-1)^k \frac{(q;q^2)_k (-q;q^2)_k^2}{(q^4;q^4)_k(-q^4;q^4)_k^2}[6k+1]q^{3k^2}
&\equiv 0\pmod{[n]},
\displaybreak[2]\\
\sum_{k=0}^{n-1}(-1)^k\frac{(q;q^2)_k (-q;q^2)_k^2}{(q^4;q^4)_k(-q^4;q^4)_k^2}[6k+1]q^{3k^2}
&\equiv 0\pmod{[n]}
\end{align*}
for all positive odd integers $n$. We plan to discuss these and other instances of such `$q$-supercongruences' in a forthcoming project.

\end{document}